\documentclass[11pt]{amsart}
\usepackage{amscd}
\usepackage{amscd}
\usepackage{latexsym}
\usepackage{pstcol,pst-plot,pst-3d}
\usepackage{mathptmx}
\usepackage{multicol}

\psset{unit=0.7cm,linewidth=0.8pt,arrowsize=2.5pt 4}

\newpsstyle{fatline}{linewidth=1.5pt}
\newpsstyle{fyp}{fillstyle=solid,fillcolor=verylight}
\definecolor{verylight}{gray}{0.97}
\definecolor{light}{gray}{0.9}
\definecolor{medium}{gray}{0.85}

\def\NZQ{\Bbb}               

\def\ZZ{{\NZQ Z}}

\def\frk{\frak}               

\def\Phi{{\frk N}}
\def\opn#1#2{\def#1{\operatorname{#2}}} 
\opn\chara{char} \opn\length{\ell} \opn\pd{pd} \opn\rk{rk}
\opn\projdim{proj\,dim} \opn\injdim{inj\,dim} \opn\rank{rank}
\opn\depth{depth} \opn\grade{grade} \opn\height{height}
\opn\embdim{emb\,dim} \opn\codim{codim}

\opn\Tr{Tr} \opn\bigrank{big\,rank}
\opn\superheight{superheight}\opn\lcm{lcm}
\opn\trdeg{tr\,deg}
\opn\reg{reg} \opn\lreg{lreg} \opn\ini{in} \opn\lpd{lpd}
\opn\size{size}\opn{\mult}{mult}
\opn\div{div} \opn\Div{Div} \opn\cl{cl} \opn\Cl{Cl}
\opn\Spec{Spec} \opn\Supp{Supp} \opn\supp{supp} \opn\Sing{Sing}
\opn\Ass{Ass} \opn\Min{Min}
\opn\Ann{Ann} \opn\Rad{Rad} \opn\Soc{Soc}
\opn\Syz{Syz} \opn\Im{Im} \opn\Ker{Ker} \opn\Coker{Coker}
\opn\Am{Am} \opn\Hom{Hom} \opn\Tor{Tor} \opn\Ext{Ext}
\opn\End{End} \opn\Aut{Aut} \opn\id{id}

\opn\nat{nat}
\opn\pff{pf}
\opn\Pf{Pf} \opn\GL{GL} \opn\SL{SL} \opn\mod{mod} \opn\ord{ord}
\opn\Gin{Gin}
\opn\Hilb{Hilb}\opn\adeg{adeg}\opn\std{std}\opn\ip{infpt}
\opn\Pol{Pol}\opn\sdepth{sdepth}\opn\sqdepth{sqdepth}\opn{\Mon}{Mon}
\opn\aff{aff} \opn\con{conv} \opn\relint{relint} \opn\st{st}
\opn\lk{lk} \opn\cn{cn} \opn\core{core} \opn\vol{vol}
\opn\link{link} \opn\star{star}
\opn\gr{gr}

\def\pot#1#2{#1[\kern-0.28ex[#2]\kern-0.28ex]}

\opn\dirlim{\underrightarrow{\lim}}
\opn\inivlim{\underleftarrow{\lim}}
\let\union=\cup
\let\sect=\cap
\let\dirsum=\oplus

\let\iso=\cong
\let\Union=\bigcup

\let\Dirsum=\bigoplus

\let\to=\rightarrow

\def\Implies{\ifmmode\Longrightarrow \else
        \unskip${}\Longrightarrow{}$\ignorespaces\fi}
\def\implies{\ifmmode\Rightarrow \else
        \unskip${}\Rightarrow{}$\ignorespaces\fi}
\def\iff{\ifmmode\Longleftrightarrow \else
        \unskip${}\Longleftrightarrow{}$\ignorespaces\fi}

\let\:=\colon
\newtheorem{Theorem}{Theorem}[section]
\newtheorem{Lemma}[Theorem]{Lemma}
\newtheorem{Corollary}[Theorem]{Corollary}
\newtheorem{Proposition}[Theorem]{Proposition}

\newtheorem{Example}[Theorem]{Example}

\let\epsilon\varepsilon
\let\phi=\varphi
\let\kappa=\varkappa
\def\qed{\ifhmode\textqed\fi
      \ifmmode\ifinner\quad\qedsymbol\else\dispqed\fi\fi}
\def\textqed{\unskip\nobreak\penalty50
       \hskip2em\hbox{}\nobreak\hfil\qedsymbol
       \parfillskip=0pt \finalhyphendemerits=0}
\def\dispqed{\rlap{\qquad\qedsymbol}}

\opn\dis{dis}
\def\pnt{{\raise0.5mm\hbox{\large\bf.}}}

\opn\Lex{Lex}

\begin{document}

\title{Stanley decompositions and partitionable simplicial complexes}

\author{J\"urgen Herzog,  Ali Soleyman Jahan and Siamak Yassemi}

\dedicatory{Dedicated to Takayuki Hibi on the occasion of his
fiftieth birthday}

\address{J\"urgen Herzog, Fachbereich Mathematik und
Informatik, Universit\"at Duisburg-Essen, Campus Essen, 45117
Essen, Germany} \email{juergen.herzog@uni-essen.de}

\address{Ali Soleyman Jahan , Fachbereich Mathematik und
Informatik, Universit\"at Duisburg-Essen, Campus Essen, 45117
Essen, Germany} \email{ali.soleyman-jahan@stud.uni-duisburg-
essen.de}

\address{Siamak Yassemi, Department of Mathematics,
University of Tehran, P.O. Box 13145448, Tehran, Iran, and
Institute for Theoretical Physics and Mathematics (IPM)}
\email{yassemi@ipm.ir}\maketitle

\begin{abstract}
We study Stanley decompositions and show that Stanley's conjecture
on Stanley decompositions implies his conjecture on partitionable
Cohen-Macaulay simplicial complexes. We also prove these conjectures
for all Cohen-Macaulay monomial ideals of codimension 2 and all
Gorenstein  monomial ideals of codimension~3.
\end{abstract}

\section*{Introduction}

In this paper we  discuss a conjecture of Stanley \cite{St2}
concerning  a combinatorial upper bound for the depth of a
$\ZZ^n$-graded module. Here we consider his conjecture only for
$S/I$, where $I$ is a monomial ideal.

Let $K$ be a field, $S=K[x_1,\ldots, x_n]$ the polynomial ring in
$n$ variables. Let $u\in S$ be a monomial and $Z$ a subset of
$\{x_1,\ldots, x_n\}$. We denote by $uK[Z]$ the $K$-subspace of
$S$ whose basis consists of all monomials $uv$ where $v$ is a
monomial in $K[Z]$. The $K$-subspace  $uK[Z]\subset S$ is called a
{\em Stanley space} of dimension $|Z|$.

Let $I\subset S$ be a monomial ideal, and denote by $I^c\subset S$
the $K$-linear subspace of $S$ spanned by all monomials which do
not belong to $I$. Then $S=I^c\dirsum I$ as a $K$-vector space,
and the residues of the monomials in $I^c$ form a $K$-basis of
$S/I$.

A decomposition $\mathcal D$ of $I^c$ as a finite direct sum of
Stanley spaces is called a {\em Stanley decomposition} of $S/I$.
Identifying $I^c$ with $S/I$, a Stanley decomposition yields a
decomposition of $S/I$ as well.  The minimal dimension of a Stanley
space in the decomposition $\mathcal D$ is called the {\em Stanley
depth} of $\mathcal D$, denoted $\sdepth ({\mathcal D})$.

We set
$
\sdepth(S/I)=\max\{\sdepth({\mathcal D})\: {\mathcal D}\; \text{is
a Stanley decomposition of $S/I$}\},
$
and call this number the {\em Stanley depth} of $S/I$.

In \cite[Conjecture 5.1]{St} Stanley conjectured the inequality
$\sdepth(S/I)\geq \depth(S/I)$. We say $I$ is a {\em Stanley
ideal}, if Stanley's conjecture holds for $I$.

Not many classes of Stanley ideals are known. Apel \cite[Corollary 3]{Ap2} showed that all monomial ideals $I$ with $\dim S/I\leq 1$ are Stanley ideals.  He also showed
\cite[Theorem 3 $\&$ Theorem 5]{Ap2} that all generic monomial
ideals and all cogeneric Cohen-Macaulay monomial ideals are Stanley
ideals, and Soleyman Jahan \cite[Proposition 2.1]{So}proved that all
monomial ideals in a polynomial ring in $n$ variables of codimension
less than or equal 1 are Stanley ideals. This implies in particular
a result of Apel which says that all monomial ideals in the
polynomial ring in three variables are Stanley ideals.

 In \cite{HePo} the authors attach to each monomial ideal a multi-complex and introduce
the concept of shellable multi-complexes. In case $I$ is a
squarefree monomial ideal, this concept of shellability coincides
with non-pure shellability introduced by Bj\"orner and Wachs
\cite{BjWa}. It is shown in \cite[Theorem 10.5]{HePo} that if $I$ is
pretty clean (see the definition in Section2), then the
multi-complex attached to $I$ is shellable and $I$ is a Stanley
ideal. The concept of pretty clean modules is a generalization of
clean modules introduced by Dress \cite{Dr}. He showed that a
simplicial complex is shellable if and only if its Stanley-Reisner
ideal is clean.

We use these results  to prove that any Cohen-Macaulay monomial
ideal of codimension 2  and that any Gorenstein monomial ideal of
codimension 3 is a Stanley ideal, see Proposition \ref{codim2} and
Theorem \ref{gorenstein}. For the proof of Proposition \ref{codim2}
we observe that the polarization of a perfect codimension 2 ideal is
shellable, and show this by using Alexander duality and result of
\cite{HeHiZh} in which it is proved that any monomial ideal with
2-linear resolution has linear quotients. The proof of Theorem
\ref{gorenstein} is based on a structure theorem for Gorenstein
monomial ideals given in \cite{BrHe1}. It also uses the result,
proved in Proposition \ref{map}, that a pretty clean monomial ideal
remains pretty clean after applying a substitution replacing the
variables by a regular sequence of monomials.

In the last section of this  paper we introduce squarefree Stanley
spaces and show in Proposition \ref{correspondance} that for a
squarefree monomial ideal $I$, the Stanley decompositions  of
$S/I$ into squarefree Stanley spaces  correspond bijectively  to
partitions into intervals of the simplicial complex whose
Stanley-Reisner ideal is the ideal $I$. Stanley calls a simplicial
complex $\Delta$ {\em partitionable} if there exists a partition
$\Delta=\Union_{i=1}^r[F_i,G_i]$  of $\Delta$ such that for all
intervals $[F_i,G_i]=\{F\in\Delta\:; F_i\subset F\subset G_i\}$
one has that $G_i$ is a facet of  $\Delta$. We show in Corollary
\ref{implies} that the Stanley-Reisner ideal $I_\Delta$ of  a
Cohen-Macaulay simplicial complex $\Delta$ is a Stanley ideal if
and only if $\Delta$ is partitionable. In other word, Stanley's
conjecture on Stanley decompositions implies his conjecture on
partitionable simplicial complexes.

\section{Stanley decompositions}

Let $S=K[x_1,\ldots,x_n]$ be a polynomial ring and $I\subset S$ a
monomial ideal.  Note that $I$ and $I^c$ as well as all Stanley
spaces are $K$-linear subspaces of $S$ with a basis which is a
subset of monomials of $S$. For any $K$-linear subspace $U\subset
S$ which is generated by monomials, we denote by $\Mon(U)$ the set
of elements in the monomial basis of $U$. It is then clear that if
$u_iK[Z_i]$, $i=1,\ldots,r$ are Stanley spaces, then
$I^c=\Dirsum_{i=1}^ru_iK[Z_i]$  if and only if $\Mon(I^c)$ is the
disjoint union of the sets $\Mon(u_iK[Z_i])$.

Usually one has infinitely many different Stanley decompositions
of $S/I$. For example if $S=K[x_1,x_2]$ and $I=(x_1x_2)$, then for
each integer $k\geq 1$ one has the Stanley decomposition
\[
\mathcal{D}_k\: \; S/I=K[x_2]\dirsum \Dirsum_{j=1}^kx_1^jK\dirsum
x_1^{k+1}K[x_1]
\]
of $S/I$. Each of these Stanley decompositions of $S/I$ has
Stanley depth $0$, while the Stanley decomposition $K[x_2]\dirsum
x_1K[x_1]$ of $S/I$ has Stanley depth 1.

Even though  $S/I$ may have infinitely many different Stanley
decompositions, all these decompositions have one property in
common, as noted in \cite[Section 2]{So}. Indeed, if $\mathcal D$
is a Stanley decomposition of $S/I$ with $s=\dim S/I$. Then the
number of Stanley sets of dimension $s$ in $\mathcal D$ is equal
to the multiplicity $e(S/I)$ of $S/I$.

There is also an upper bound for $\sdepth(S/I)$ known, namely
\[
\sdepth(S/I)\leq \min\{\dim S/P \: P\in \Ass(S/I)\}.
\]
see \cite[Section 3]{Ap2}. Note that for  $\depth(S/I)$  the same
upper bound is valid. As a consequence of these observations one
has

\begin{Corollary} \label{cm} Let $I\subset S$ be a monomial ideal such that
$S/I$ is Cohen-Macaulay. Then the following conditions are
equivalent:
\begin{enumerate}
\item[(a)]  $I$ is a Stanley ideal.

\item[(b)] There exists a Stanley decomposition $\mathcal D$ of
$S/I$ such that each Stanley space in $\mathcal D$ has dimension
$d=\dim S/I$.

\item[(c)] There exists a Stanley decomposition $\mathcal D$ of
$S/I$ which has $e(S/I)$ summands.
\end{enumerate}
\end{Corollary}

The following result will be needed later in Section 2.

\begin{Proposition}\label{complete} Let  $I\subset S$ be  a monomial complete intersection ideal.
Then  $S/I$ is clean. In particular, $I$ is a Stanley  ideal.
\end{Proposition}

\begin{proof} Let $u\in S$ be a monomial. We call $\supp(u)=\{x_i\: x_i\;
\text{divides}\; u\}$ the {\em support} of $u$. Now let
$G(I)=\{u_1,\ldots,u_m\}$ be the unique minimal set of monomial
generators of $I$. By our assumption, $u_1,\ldots, u_m$ is a regular
sequence. This implies that $\supp(u_i)\sect\supp(u_j)=\emptyset$
for all $i\neq j$.

It follows from the definition of the polarization of a monomial
ideal (see for example \cite{So}), that for the polarized ideal
$I^p=(u_1^p., \ldots,u_m^p)$ one again has
$\supp(u_i^p)\sect\supp(u_j^p)=\emptyset$ for all $i\neq j$.

Thus $J=I^p$ is a squarefree monomial ideal generated by the
regular sequence of monomials $v_1,\ldots, v_m$ with $v_i=u_i^p$
for all $i$.

Let $\Delta$ be the simplicial complex whose Stanley-Reisner ideal
$I_\Delta$ is equal to $J$. The {\em Alexander dual} $\Delta^\vee$
of $\Delta$ is defined to be the simplicial complex whose faces
are $\{[n]\setminus F\:\; F\not\in \Delta\}$. The Stanley-Reisner
ideal of  ${\Delta^\vee}$ is minimally generated by all monomials
$x_{i_1}\cdots x_{i_k}$ where $(x_{i_1},\ldots, x_{i_k})$ is a
minimal prime ideal of $I_\Delta$.

In our case it follows that $I_{\Delta^\vee}$ is minimally generated
by the monomial of the form $x_{i_1}\ldots x_{i_m}$ where
$x_{i_j}\in\supp(v_j)$ for $j=1,\ldots,m$. Thus we see that
$I_{\Delta^\vee}$ is the matriodal ideal of the transversal matroid
attached to the sets $\supp(v_1),\ldots, \supp(v_m)$, see
\cite[Section 5]{CoHe}. In \cite[Lemma~1.3]{HeTa} and \cite[Section
5]{CoHe} it is shown that any polymatroidal ideal has linear
quotients, and this implies that $\Delta$ is a shellable simplicial
complex, see for example \cite[Theorem 1.4]{HeHiZh1}. Hence by the
theorem of Dress quoted in the next section, $S/I_\Delta$ is clean.
Now we use the result in \cite[Theorem 3.10]{So} which says that a
monomial ideal is pretty clean (see the definition in Section 2) if
and only if its polarization is clean. Therefore we conclude that
$S/I$ is pretty clean. Since all prime ideals in a pretty clean
filtration are associated prime ideals of $S/I$ (see \cite[Corollary
3.4]{HePo}) and since $S/I$ is Cohen-Macaulay, the prime ideals in
the filtration are minimal. Hence $S/I$ is clean. Thus we conclude
from \cite[Theorem 6.5]{HePo} that $I$ is Stanley ideal.
\end{proof}

\begin{Corollary}
\label{easy} Let $I\subset S$ be a monomial ideal with $\depth
S/I\geq n-1$. Then $I$ is a Stanley ideal.
\end{Corollary}

\begin{proof}
The assumption implies that $I$ is a principal ideal. Thus the
assertion follows from Proposition \ref{complete}.
\end{proof}

With the same techniques as in the proof of Proposition
\ref{complete} we can show

\begin{Proposition}
\label{codim2} Let $I\subset S$ be a monomial  ideal which is
perfect and of codimension 2. Then $S/I$ is clean. In particular,
$I$ is a Stanley ideal.
\end{Proposition}

\begin{proof} We will show that the polarized ideal $I^p$ defines
a shellable simplicial complex. Then, as in the proof of Proposition
\ref{complete}, it follows that $S/I$ is clean. Note that $I^p$ is a
perfect squarefree monomial ideal of codimension 2. Let $\Delta$ be
the simplicial complex defined by $I^p$. By the Eagon--Reiner
theorem \cite{ER} and a result of Terai \cite{T}, the ideal
$I_{\Delta^\vee}$ has a 2-linear resolution. Now we use the fact,
proved in \cite[Theorem 3.2]{HeHiZh}, that an ideal with 2-linear
resolution has linear quotients which in turn implies that $\Delta$
is shellable, as desired. \end{proof}

Combining the preceding results with Apel's result according to
which   all monomial ideals with $\dim S/I\leq 1$ are Stanley
ideals we obtain

\begin{Corollary} Let $I\subset S$ be a monomial ideal. If $n\leq
4$ and  $S/I$ is Cohen-Macaulay, then $I$ is a
Stanley ideal.
\end{Corollary}

\section{Gorenstein monomial ideals of codimension 3}

As the main result of  this section we will show

\begin{Theorem}
\label{gorenstein} Each  Gorenstein monomial ideal of codimension
3 is a Stanley ideal. \end{Theorem}

 The proof of this result is based on the following structure
 theorem that can be found in \cite{BrHe1}.

\begin{Theorem}
\label{structure} Let $I\subset S$ be a monomial Gorenstein ideal
of codimension 3. Then $G(I)$ is an odd number, say $|G(I)|=2m+1$,
and there exists a regular sequence of monomials $u_1,\ldots
u_{2m+1}$ in $S$ such that  $$G(I)= \{u_iu_{i+1}\cdots
u_{i+m-1}\:\; i=1,\ldots, 2m+1\},$$ where $u_i=u_{i-2m-1}$
whenever  $i>2m+1$.
 \end{Theorem}

In order to apply this theorem we need another result. Let
$I\subset S$ be monomial ideal. According to \cite{HePo}, $S/I$ is
called {\em pretty clean}, if there exists a chain of monomial
ideals  such that
\begin{enumerate}
\item[(a)] for all $j$ one has $I_j/I_{j+1}\iso S/P_j$ where $P_j$
is a monomial prime ideal;

\item[(b)] for all $i<j$ such that $P_i\subset P_j$, it follows
that $P_i=P_j$.
\end{enumerate}

Dress \cite{Dr} calls the ring $S/I$  {\em clean}, if there exists a
chain of ideals as above such that all the $P_i$ are minimal prime
ideals of $I$. By an abuse of notation we call $I$ (pretty) clean if
$S/I$ is (pretty) clean. Obviously, any clean ideal is pretty clean.
In \cite[Theorem 6.5]{HePo} it is shown that if $I$ is pretty clean,
then $I$ is a Stanley ideal, while Dress showed \cite[Section 4]{Dr}
that if $I=I_\Delta$ for some simplicial complex $\Delta$, then
$\Delta$ is shellable if and only if $I_\Delta$ is clean. In
particular, it follows that $I_\Delta$ is a Stanley ideal, if
$\Delta$ is shellable.

We now show

\begin{Proposition}
\label{map} Let $I\subset T= K[y_1,\ldots, y_r]$ be a monomial
ideal such that $T/I$ is (pretty) clean. Let $u_1,\ldots, u_r\in
S=K[x_1,\ldots, x_n]$ be a regular sequence of monomials, and  let
$\varphi \:\; T\to S$ be the $K$-algebra homomorphism with
$\varphi(y_j)=u_j$ for $j=1,\ldots, r$. Then  $S/\varphi(I)S$
 is (pretty) clean.
\end{Proposition}

\begin{proof}
Let $I=I_0\subset I_1\subset\cdots \subset I_m=T$ be a pretty clean
filtration $\mathcal F$ of $T/I$ with $I_k/I_{k+1}=T/P_k$ for all
$k$.

Observe that the $K$-algebra homomorphism $\varphi\: T\to S$ is
flat, since $u_1,\ldots, u_r$ is a regular sequence. Hence if we
set $J_k=\phi(I_k)S$ for $k=1,\ldots,m$, then we obtain the
filtration $\varphi(I)S=J_0\subset J_1\subset \cdots \subset
J_m=S$ with $J_k/J_{k+1}\iso S/\varphi(P_k)S$.

Suppose $P_k=(y_{i_1},\ldots, y_{i_k})$, then
$\varphi(P_k)S=(u_{i_1},\ldots, u_{i_k})$. In other words,
$\varphi(P_k)S$ is a monomial complete intersection, and hence by
Proposition \ref{complete} we have that $S/\varphi(P_k)S$ is clean.
Therefore there exists a prime filtration $J_k=J_{k_0}\subset
J_{k_1}\subset\cdots \subset J_{k_{r_k}}=J_{k+1}$ such that
$J_{k_i}/J_{k_{i+1}}\iso S/P_{k_i}$ where $P_{k_i}$ is a minimal
prime ideal of $\varphi(P_k)S$. Since
$\varphi(P_k)S=(u_{i_1},\ldots, u_{i_{t_k}})S$ is a complete
intersection, all minimal prime ideals of $\varphi(P_k)$ have height
$t_k$.

Composing the prime filtrations of the $J_k/J_{k+1}$, we obtain a
prime filtration of $S/\varphi(I)S$. We claim that this prime
filtration is (pretty) clean.   In fact, let $P_{k_i}$ and
$P_{\ell_j}$ be two prime ideals in the support of this filtration.
We have to show: if  $P_{k_i}\subset P_{\ell_j}$ for $k<\ell$, or
$P_{k_i}\subset P_{\ell_j}$ for $k=\ell$ and $i<j$, then $P_{k_i}=
P_{\ell_j}$. In  case $k=\ell$, we have
$\height(P_{k_i})=\height(P_{\ell_j})=t_k$, and the assertion
follows. In case $k<\ell$, by using the fact that $\mathcal F$ is a
pretty clean filtration, we have that  $P_k=P_\ell$ or $P_k\not
\subset P_\ell$. In the first case, the prime ideals $P_{k_i}$ and
$P_{\ell_j}$ have the same height, and the assertion follows. In the
second case there exists a variable $y_g\in P_k\setminus P_\ell$.
Then the monomial $u_g$ belongs to $\varphi(P_k)S$ but not to
$\varphi(P_\ell)S$. This implies that $P_{k_i}$ contains a variable
which belongs to the support of $u_g$. However this variable cannot
be a generator of $P_{\ell_j}$, because the support of $u_g$ is
disjoint of the support of all the monomial generators of
$\varphi(P_\ell)S$. This shows that $P_{k_i}\not\subset P_{\ell_j}$.
\end{proof}

\begin{Corollary}
\label{transfer} Let $\Delta$ be a shellable simplicial complex and
$I_\Delta\subset T=K[y_1,\ldots,y_r]$ its Stanley-Reisner ideal.
Furthermore, let $u_1,\ldots,u_r\subset S=K[x_1,\ldots,x_n]$ be a
regular sequence of monomials, and let $\varphi(y_i)=u_i$ for
$i=1,\ldots,r$. Then $\varphi(I_\Delta)S$ is a Stanley ideal.
\end{Corollary}

\begin{proof} By the theorem of Dress, the ring $T/I_{\Delta}$ is
clean. Therefore, $S/\varphi(I_\Delta)S$ is again clean, by
Proposition \ref{map}. In particular, $S/\varphi(I_\Delta)S$  is
pretty clean which according to \cite[Theorem 6.5]{HePo}  implies
that $\varphi(I_\Delta)S$ is a Stanley ideal.
\end{proof}

\begin{proof}[Proof of Theorem \ref{gorenstein}] Let $\Delta$ be
the simplicial complex whose Stanley-Reisner ideal
$$I_\Delta\subset T=K[y_1,\ldots, y_{2m+1}]$$  is generated by the
monomials $y_iy_{i+1}\cdots y_{i+m-1}$, $i=1,\ldots, 2m+1$, where
$y_i=y_{i-2m-1}$ whenever  $i>2m+1$, and let $u_1,\ldots,
u_{2m+1}\subset S=K[x_1,\ldots,x_n]$ be the regular sequence given
in Theorem \ref{gorenstein}. Then we have $I=\varphi(I_\Delta)S$
where $\varphi(y_j)=u_j$ for all $j$. Therefore, by Corollary
\ref{transfer}, it suffices to show that $\Delta$ is shellable.

Identifying the vertex set of $\Delta$ with $[2m+1]=\{1,\ldots,
2m+1\}$ and observing that $I_\Delta$ is of codimension 3, it is
easy to see that $F\subset [2m+1]$ is a facet of $\Delta$ if and
only if $F=[2m+1]\setminus\{a_1,a_2,a_3\}$ with
\[
a_2-a_1< m+1,\qquad a_3-a_2<m+1,\qquad a_3-a_1>m.
\]
We denote the facet $[2m+1]\setminus\{a_1,a_2,a_3\}$ by
$\mathrm{F}(a_1,a_2,a_3)$

We will show that $\Delta$ is shellable with respect to the
lexicographic order. Note that
$\mathrm{F}(a_1,a_2,a_3)<\mathrm{F}(b_1,b_2,b_3)$ in the
lexicographic order, if and only if either $b_1<a_1$, or $b_1=a_1$
and $b_2<a_2$, or $a_1=b_1$, $a_2=b_2$ and $b_3<a_3$.

In order to prove that $\Delta$ is shellable we have to show: if
$F=\mathrm{F}(a_1,a_2,a_3)$ and $G=\mathrm{F}(b_1,b_2,b_3)$ with
$F<G$,  then  there exists $c\in G\setminus F$ and some facet $H$
such that $H<G$ and $G\setminus H=\{c\}$.

We know that $|G\setminus F|\le 3$. If $|G\setminus F|=1$, then
there is nothing to prove. In the following we discuss the cases
$|G\setminus F|=2$ and $|G\setminus F|=3$. The discussion of these
cases is somewhat tedious but elementary. For the convenience of
the reader we list all the possible cases.

\medskip
\noindent Case 1: $|G\setminus F|=2$.

\medskip
 (i) If $b_1=a_1<b_2<a_2$, then we choose
$H=(G\setminus\{a_2\})\cup\{b_2\}$.

 (ii) If $b_1<b_2=a_1$ or $b_1<b_2<a_1<a_2=b_3<a_3$, then we choose
$H=(G\setminus\{a_3\})\cup\{b_1\}$.

(iii) If $b_1<a_1<b_2<a_2=b_3<a_3$, we consider the following two
subcases:

\begin{verse}

for $a_3-b_2<m+1$, we choose $H=(G\setminus\{a_3\})\cup\{b_3\}$.

for $a_3-b_2\ge m+1$,  we choose
$H=(G\setminus\{a_3\})\cup\{b_1\}$.

\end{verse}

(iv) If $b_1<a_1<a_2=b_2<b_3<a_3$, then we choose
$H=(G\setminus\{a_3\})\cup\{b_3\}$.

(v) If $b_1<a_1<a_2=b_2<a_3<b_3$ or $b_1<a_1<a_2<a_3=b_2<b_3$,
then we choose $H=(G\setminus\{a_1\})\cup\{b_1\}$.

\medskip
\noindent Case 2: $|G\setminus F|=3$.
\medskip

(i) If $b_1<a_1<a_2<a_3<b_3$, then we choose
$H=(G\setminus\{a_1\})\cup\{b_1\}$.

(ii) If $b_1<b_2<b_3<a_1<a_2<a_3$ or $b_1<b_2<a_1<a_2<a_3$ and
$a_1<b_3$, then we choose $H=(G\setminus\{a_1\})\cup\{b_2\}$.

(iii) If $b_1<a_1<b_2<b_3<a_2<a_3$, then we choose
$H=(G\setminus\{a_2\})\cup\{b_3\}$.

(iv) If $b_1<a_1<b_2<a_2<b_3<a_3$,  we consider the following two
subcases:

\begin{verse}

for $a_3-b_2<m+1$, we choose $H=(G\setminus\{a_3\})\cup\{b_3\}$.

for $a_3-b_2\ge m+1$,  we choose
$H=(G\setminus\{a_3\})\cup\{b_1\}$.

\end{verse}

(iv) If $b_1<a_1<a_2<b_2<b_3<a_3$, then we choose
$H=(G\setminus\{a_3\})\cup\{b_3\}$.
\end{proof}

Combining the result of Theorem \ref{gorenstein} with the result of Apel \cite[Corollary 3]{Ap2} we obtain

\begin{Corollary}
\label{need} Let $I\subset S$ be monomial ideal. If $n\leq 5$ and
$S/I$ is Gorenstein, then  $I$ is a Stanley ideal.
\end{Corollary}

\section{Squarefree Stanley decompositions and partitions of simplicial complexes}

 A Stanley space $uK[Z]$ is called a {\em
squarefree Stanley space}, if $u$ is a squarefree monomial and
$\supp(u)\subseteq Z$. We shall use the following notation: for
$F\subseteq [n]$ we set $x_F=\prod_{i\in F}x_i$ and $Z_F=\{x_i\:\,
i\in F\}$. Then a Stanley space is squarefree if and only if it is
of the form $x_FK[Z_G]$ with $F\subseteq G\subseteq [n]$.

A Stanley decomposition of $S/I$ is  called a {\em squarefree
Stanley decomposition} of $S/I$, if all Stanley spaces in the
decomposition are squarefree.

\begin{Lemma}
\label{squarefree} Let $I\subset S$ be a monomial ideal. The
following conditions are equivalent:
\begin{enumerate}
\item[(a)] $I$ is a squarefree monomial ideal.

\item[(b)] $S/I$ has a squarefree Stanley decomposition.
\end{enumerate}
\end{Lemma}

\begin{proof} (a)\implies (b): We may view $I$ as the
Stanley-Reisner ideal of some simplicial complex $\Delta$.  With
each $F\in \Delta$ we associate the squarefree Stanley space
$x_FK[Z_F]$. We claim that $\Dirsum_{F\in\Delta}x_FK[Z_F]$ is a
(squarefree) Stanley decomposition of $S/I$. Indeed, a monomial
$u\in S$ belongs to $I^c$ if and only if $\supp(u)\in\Delta$, and
these monomial form a $K$-basis for $I^c$. On the other hand, a
monomial $u\in S$ belongs to $x_FK[Z_F]$ if and only if
$\supp(u)=F$. This shows that $I^c=\Dirsum_{F\in\Delta}x_FK[Z_F]$.

(b)\implies (a): Let $\Dirsum_i u_iK[Z_i]$ be a squarefree Stanley
decomposition of $S/I$. Assume that $I$ is not a squarefree monomial
ideal. Then there exists $u\in G(I)$ which is not squarefree and we
may assume that $x_1^2|u$. Then $u'=u/x_1\in I^c$, and hence there
exists $i$ such that $u'\in u_iK[Z_i]$. Since $x_1|u'$ it follows
that  $x_1\in Z_i$. Therefore $u\in u_iK[Z_i]\subset I^c$, a
contradiction.
\end{proof}

Let $\Delta$ be a simplicial complex of dimension $d-1$ on the
vertex set $V=\{x_1,\ldots,x_n\}$. A subset ${\mathcal I}\subset
\Delta$ is called an {\em interval}, if there exits faces $F,G\in
\Delta$ such that ${\mathcal I}=\{H\in \Delta \: F\subseteq
H\subseteq G\}$. We denote this interval given by $F$ and $G$ also
by $[F,G]$ and call $\dim G-\dim F$ the {\em rank} of the interval.
A {\em partition} $\mathcal P$ of $\Delta$ is a presentation of
$\Delta$ as a disjoint union of intervals. The $r$-vector of
$\mathcal P$ is the integer vector $r=(r_0, r_1,\ldots, r_{d})$
where $r_i$ is the number of intervals of rank $i$.

\begin{Proposition}
\label{correspondance} Let ${\mathcal P}\:\;
\Delta=\Union_{i=1}^r[F_i,G_i]$ be a partition of $\Delta$. Then
\begin{enumerate}
\item[(a)]  $D({\mathcal P})=\Dirsum_{i=1}^rx_{F_i}K[Z_{G_i}]$ is
squarefree Stanley decomposition of $S/I$.

\item[(b)] The map ${\mathcal P}\mapsto D({\mathcal P})$
establishes a bijection between  partitions of $\Delta$ and
squarefree Stanley decompositions of $S/I$.
\end{enumerate}
\end{Proposition}

\begin{proof} (a) Since each $x_{F_i}K[Z_{G_i}]$ is a squarefree
Stanley space  it suffices to show that $I^c$ is indeed the direct
sum of the Stanley spaces $x_{F_i}K[Z_{G_i}]$. Let $u\in \Mon(I^c)$;
then $H=\supp(u)\in \Delta$. Since $\mathcal P$ is a partition of
$\Delta$ it follows that $H\in [F_i,G_i]$ for some $i$. Therefore,
$u=x_{F_i}u'$ for some monomial $u'\in K[Z_{G_i}]$. This implies
that $u\in x_{F_i}K[Z_{G_i}]$. This shows that $\Mon(I^c)$ is the
union of sets $\Mon(x_{F_i}K[Z_{G_i}])$. Suppose there exists a
monomial $u\in x_{F_i}K[Z_{G_i}]\sect x_{F_j}K[Z_{G_j}]$. Then
$\supp(u)\in [F_i,G_i]\sect [F_j,G_j]$. This is only possible if
$i=j$, since $\mathcal P$ is partition of $\Delta$.

(b) Let  $[F_i,G_i]$  and $[F_j,G_j]$ be two intervals. Then
$x_{F_i}K[Z_{G_i}]=x_{F_j}K[Z_{G_j}]$ if and only if
$[F_i,G_i]=[F_j,G_j]$. Indeed, if
$x_{F_i}K[Z_{G_i}]=x_{F_j}K[Z_{G_j}]$, then $x_{F_j}\in
x_{F_i}K[Z_{G_i}]$, and hence $x_{F_i}|x_{F_j}$. By symmetry we
also have $x_{F_j}|x_{F_i}$. In other words, $F_i=F_j$, and it
also follows that $K[Z_{G_i}]=K[Z_{G_j}]$. This implies $G_i=G_j$.
These considerations show that ${\mathcal P}\mapsto D({\mathcal
P})$ is injective.

On the other hand, let ${\mathcal D}\:\;
S/I=\Dirsum_{i=1}^rx_{F_i}K[Z_{G_i}]$ be an arbitrary squarefree
Stanley decomposition of $S/I$. By the definition of a squarefree
Stanley set we have $F_i\subseteq G_i$, and since
$x_{F_i}K[Z_{G_i}]\subset I^c$, it follows that $G_i\in \Delta$.
Hence $[F_i,G_i]$ is an interval of $\Delta$, and a squarefree
monomial $x_F$ belongs to $x_{F_i}K[Z_{G_i}]$ if and only if $F\in
[F_i,G_i]$.

Let $F\subset \Delta$ be an arbitrary face. Then $x_F\in
\Mon(I^c)=\Union_{i=1}^r\Mon(x_{F_i}K[Z_{G_i}])$. Hence the
squarefree monomial $x_F$ belongs to $x_{F_i}K[Z_{G_i}]$ for some
$i$, and hence $F\in [F_i,G_i]$. This shows that
$\Union_{i=1}^r[F_i,G_i]=\Delta$. Suppose $F\in [F_i,G_i]\sect
[F_j,G_j]$. Then $x_F\in x_{F_i}K[Z_{G_i}]\sect x_{F_j}K[Z_{G_j}]$,
a contradiction. Hence we see that ${\mathcal P}\:\;
\Delta=\Union_{i=1}^r[F_i,G_i]$ is a partition of $\Delta$ with
$D({\mathcal P})=\mathcal D$.
\end{proof}

Now let $I\subset S$ be a squarefree monomial ideal. Then we set
\[ \sqdepth(S/I)=\max\{\sdepth({\mathcal D})\: {\mathcal D}\; \text{is
a squarefree Stanley decomposition of $S/I$}\},
\]
and call this number the {\em squarefree Stanley depth} of $S/I$.

As the main result of this section we have

\begin{Theorem}
\label{ss=s} Let $I\subset S$ be a squarefree monomial ideal. Then
$\sqdepth(S/I)=\sdepth(S/I)$.
\end{Theorem}

\begin{proof}
Let $\mathcal D$ be any Stanley decomposition of $S/I$, and let
$\Delta$ be the simplicial complex with $I=I_\Delta$. For each
$F\in \Delta$ we have $x_F\in I^c$. Hence there exists a summand
$uK[Z]$ with $x_F\in uK[Z]$. Since $x_F$ is squarefree it follows
that $u=x_G$ is squarefree and $F\subseteq G\union Z$. Let
${\mathcal D}'$ the sum of those Stanley spaces $uK[Z]$ in
$\mathcal D$ for which $u$ is a squarefree monomial. Then this sum
is direct. Therefore the intervals $[G,G\union Z]$ corresponding
to the summands in ${\mathcal D}'$ are pairwise disjoint. On the
other hand these intervals cover $\Delta$, as we have seen before,
and hence form a partition of $\mathcal P$ of $\Delta$. It follows
from the construction of $\mathcal P$ that $\sqdepth D({\mathcal
P})\geq \sdepth \mathcal D$. This shows that $\sqdepth(S/I)\geq
\sdepth(S/I)$. The other inequality $\sqdepth(S/I)\leq
\sdepth(S/I)$ is obvious.
\end{proof}

\begin{Corollary}
\label{partition} Let $\Delta$ be a simplicial complex. Then the
following conditions are equivalent:
\begin{enumerate}
\item[(a)] $I_\Delta$ is a Stanley ideal.

\item[(b)] There exists a partition
$\Delta=\Union_{i=1}^r[F_i,G_i]$ with $|G_i|\geq \depth K[\Delta]$
for all $i$.
\end{enumerate}
\end{Corollary}

Let $\Delta$ be a simplicial complex and ${\mathcal F}(\Delta)$
its set of facets. Stanley calls a simplicial complex $\Delta$
{\em partitionable} if there exists a partition
$\Delta=\Union_{i=1}^r[F_i,G_i]$ with ${\mathcal
F}(\Delta)=\{G_1,\ldots, G_r\}$. We call a partition with this
property a {\em nice partition}.  Stanley  conjectures
\cite[Conjecture 2.7]{St1} (see also \cite[Problem 6]{St2}) that
each Cohen-Macaulay simplicial complex is partitionable. In view
of Corollary \ref{cm} it follows that the conjecture of Stanley
decompositions implies the conjecture on partitionable simplicial
complexes. More precisely we have

\begin{Corollary}
\label{implies} Let $\Delta$ be a Cohen-Macaulay simplicial
complex with $h$-vector $(h_0,h_1,\ldots,h_{d})$. Then the
following conditions are equivalent:
\begin{enumerate}
\item[(a)] $I_\Delta$ is a Stanley ideal.

\item[(b)] $\Delta$ is partitionable.

\item[(c)] $\Delta$ admits a partition whose $r$-vector satisfies
$r_i=h_{d-i}$ for $i=0,\ldots,d$.

\item[(d)] $\Delta$ admits a partition into $e(K[\Delta])$
intervals.
\end{enumerate}
Moreover, any nice partition of $\Delta$ satisfies the conditions
{\em (c)} and {\em (d)}.
\end{Corollary}

\begin{proof} (a)\iff (b) follows from Corollary \ref{partition}.
In order to prove the implication (b)\implies (c), consider a nice
partition $\Delta=\Union_{i=1}^r[F_i,G_i]$ of $\Delta$. From this
decomposition the $f$-vector of $\Delta$ can be computed by the
following formula
\[
\sum_{i=0}^df_{i-1}t^i=\sum_{i=0}^dr_it^{d-i}(1+t)^i.
\]
On the other hand one has
 \[
\sum_{i=0}^df_{i-1}t^i=\sum_{i=0}^dh_it^{i}(1+t)^{d-i},
\]
see \cite[p.\ 213]{BrHe}. Comparing coefficients the assertion
follows.

The implication (c)\implies (d) follows from the fact that
$e(K[\Delta])=\sum_{i=0}^dh_i$, see
\cite[Proposition~4.1.9]{BrHe}. Finally (d) \implies (a) follows
from Corollary \ref{cm}.
\end{proof}

We conclude this section with some explicit examples. Recall that
constructibility,  a generalization of shellability, is  defined
recursively as follows: (i) a simplex is constructible, (ii) if
$\Delta_1$ and $\Delta_2$ are $d$-dimensional constructible
complexes and their intersection is a $(d-1)$-dimensional
constructible complex, then their union is constructible. In this
definition, if in the recursion we restrict  $\Delta_2$ always to be
a simplex, then the definition becomes equivalent to that of (pure)
shellability. The notion of constructibility for simplicial
complexes  appears in \cite{St3}. It is known and easy to see that
\begin{center}
Shellable $\Rightarrow$ constructible $\Rightarrow$ Cohen-Macaulay.
\end{center}

Since any shellable simplicial complex is partitionable (see
\cite[p.\ 79]{St1}), it is natural to ask whether any constructible
complex is partitionable?  This question is a special case of
Stanley's conjecture that says that Cohen-Macaulay simplicial
complexes are partitionable. We do not know the answer yet! In the
following we present some examples where the complexes are not
shellable or are not Cohen-Macaulay but the ideals related to these
simplicial complexes are Stanley ideals.

\begin{Example}{\em The following example of a simplicial complex  is due to Masahiro Hachimori
\cite{Ha}. The simplicial complex $\Delta$ described by the next
figure is 2-dimensional, non shellable but constructible. It is
constructible, because if we divide the simplicial complex by the
bold line, we obtain two shellable complexes, and their intersection
is a shellable 1-dimensional simplicial complex.

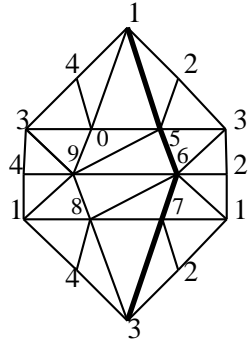
\psset{unit=0.15cm}
\begin{pspicture}(20,-10)(-35,23)

\put(1.1,10){$3$}

\put(8.3,8.5){\footnotesize$0$}

\put(14.7,8.5){\footnotesize$5$}

\put(20.4,10){$3$}

\put(0.5,6){$4$}

\put(5.7,7.2){\footnotesize$9$}

\put(15.4,7){\footnotesize$6$}

\put(20.4,6){$2$}

\put(0.5,2){$1$}

\put(6,2.5){\footnotesize$8$}

\put(15,2.5){\footnotesize$7$}

\put(20.4,2){$1$}

\put(11,19.5){$1$}

\put(11,-8.5){$3$}

\put(5.5,15){$4$}

\put(16,15){$2$}

\put(5.5,-4){$4$}

\put(16,-3.5){$2$}



 \psline(2,10)(7.8,10)\psline(2,10)(6.5,14.5) \psline(2,10)(1.8,6)
 \psline(2,10)(6.2,6)

\psline(6.5,14.5)(7.8,10)\psline(11,19)(7.8,10)\psline(13.9,10)(7.8,10)
\psline(6.2,6)(7.8,10)

\psline(6.2,6)(2,10)\psline(6.2,6)(1.8,6)\psline(6.2,6)(13.9,10)\psline(6.2,6)(15.4,6)
\psline(6.2,6)(7.7,2)\psline(6.2,6)(1.8,2)

\psline(6.5,14.5)(11,19)\psline(1.8,6)(1.8,2)

\psline(7.7,2)(1.8,2)\psline(7.7,2)(6.5,-2.5)\psline(7.7,2)(11,-7)\psline(7.7,2)(15.4,6)

\psline(13.8,2)(19.8,2)\psline(14.1,2)(15.5,-2.5)
\psline(7.7,2)(13.8,2)

\psline(15.4,6)(19.8,2)\psline(15.4,6)(19.8,6)\psline(15.4,6)(19.7,10)
\psline[linewidth=2pt](15.4,6)(13.9,10)

\psline[linewidth=2pt](13.9,10)(11,19)\psline(15.5,14.5)(13.9,10)\psline(13.9,10)(19.7,10)
\psline(11,19)(15.5,14.5)\psline(15.5,14.5)(19.7,10)\psline(19.8,6)(19.7,10)
\psline(19.8,6)(19.8,2)
\psline(15.5,-2.5)(19.8,2)\psline(15.5,-2.5)(11,-7)\psline(6.5,-2.5)(11,-7)
\psline(6.5,-2.5)(1.8,2)

\psline[linewidth=2pt](15.4,6)(11,-7)

\end{pspicture}

\noindent
 Indeed we can write $\Delta=\Delta_1\union \Delta_2$ where
the shelling order of the facets of $\Delta_1$ is given by: $$148,
149, 140, 150, 189,  348, 349, 378, 340, 390,  590, 569, 689,678, $$
and that of $\Delta_2$ is given by: $$125,126, 127, 167,
235,236,237, 356.$$

\noindent We use the following principle to construct a partition of
$\Delta$: suppose that $\Delta_1$ and $\Delta_2$ are $d$-dimensional
partitionable simplicial complexes, and that $\Gamma=\Delta_1\sect
\Delta_2$ is $(d-1)$-dimensional pure simplicial complex. Let
$\Delta_1=\Union_{i=1}^r[K_i,L_i]$ be a nice partition of
$\Delta_1$, and $\Delta_2=\Union_{i=1}^s[F_i,G_i]$ a nice partition
of $\Delta_2$. Suppose that for each $i$, the set
$[F_i,G_i]\setminus \Gamma$ has a unique minimal element $H_i$. Then
$\Delta_1\union \Delta_2=\Union_{i=1}^r[K_i,L_i]\union
\Union_{i=1}^s[H_i,G_i]$ is a nice partition of $\Delta_1\union
\Delta_2$. Notice that $[F_i,G_i]\setminus \Gamma$ has a unique
minimal element if and only if for all $F\in [F_i,G_i]\sect \Gamma$
there exists a facet $G$ of $\Gamma$ with $F\subseteq G\subset G_i$.

Suppose  that $\Delta_2$ is shellable with shelling $G_1,\ldots,
G_s$. Let $F_i$ be  the unique minimal subface of $G_i$ which is
not a subface of any $G_j$ with $j<i$. Then
$\Delta_2=\Union_{i=1}^s[F_i,G_i]$ is the nice partition induced
by this shelling. The above discussions then show that
$\Delta_1\union \Delta_2$ is partitionable, if for all $i$ and all
$F\in \Gamma$ such that $F\subset G_i$ and $F\not \subset G_j$ for
$j<i$, there exists a facet $G\in\Gamma$ with $F\subseteq G\subset
G_i$.

\medskip
In our particular case the shelling of $\Delta_1$ induces the
following partition of $\Delta_1$: \[[\emptyset, 148], [9,149],
[0,140], [5,150], [89,189], [3,348], [39,349],[7,378], \]
\[[30,340],[90,390],[59,590],[6,569],[68,689],[67,678],\]
and  the shelling of $\Delta_2$ induces the following partition of
$\Delta_2$:
\[
[\emptyset, 125], [6,126],
[7,127],[67,167],[3,235],[36,236],[37,237],[56,356].
\]
The facets of  $\Gamma=\Delta_1\sect \Delta_2$ are:  $15,56,67,73.$

The restriction of the intervals of this partition  of $\Delta_2$
to the complement of $\Gamma$ do not all give intervals. For
example we have $[6,126]\setminus \Gamma=\{16,26,126\}$.  This set
has two minimal elements, and hence is not an interval.  On the
other hand, the following partition of $\Delta_2$ (which is not
induced from a shelling)
\[
[\emptyset,237],[1,125],[5,356],[6,167],[17,127],[25,235],[26,126],[36,236]
\]
restricted to the complement of $\Gamma$ yields the following
intervals
\[
[2,237],
[12,125],[35,356],[16,167],[17,127],[25,235],[26,126],[36,236],
\]
which together with the intervals of the partition of $\Delta_1$
give us a partition of $\Delta$.

}
\end{Example}

\begin{Example}{\em (The Dunce hat)
The Dunce hat is the topological space obtained from the solid
triangle $abc$ by identifying the oriented edges $\vec{ab}$,
$\vec{bc}$ and $\vec{ac}$. The following is a triangulation of the
Dunce hat using 8 vertices.

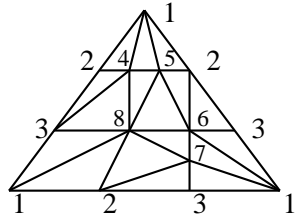
\psset{unit=0.2cm}
\begin{pspicture}(-13,-1)(-20,17)





\put(2,0.5){$1$} \put(8,0.5){$2$} \put(14,0,5){$3$}
\put(20,0.5){$1$} \put(3.5,5.5){$3$} \put(6.5,10){$2$}
\put(17.9,5.5){$3$} \put(14.9,10){$2$} \put(12,13.2){$1$}

\put(8.7,6.3){\footnotesize 8} \put(14.3,6.3){\footnotesize 6}
\put(14.1,4){\footnotesize 7} \put(9,10.3){\footnotesize 4}
\put(12.1,10.3){\footnotesize 5}

\psline(1.77,2)(10.77,14)\psline(10.77,14)(19.77,2)\psline(19.77,2)(1.77,2)

\psline(4.77,6)(16.77,6)\psline(7.77,10)(13.77,10)

\psline(7.77,2)(11.77,10)\psline(13.77,2)(13.77,10)

\psline(10.77,14)(9.77,10)\psline(10.77,14)(11.77,10)

\psline(11.77,10)(13.77,6)

\psline(9.77,10)(9.77,6) \psline(9.77,10)(4.77,6)

\psline(9.77,6)(1.77,2)\psline(9.77,6)(13.77,4)

\psline(13.77,4)(7.77,2)\psline(13.77,4)(19.77,2)
\psline(13.77,6)(19.77,2)










\end{pspicture}

\noindent The facets  arising from this triangulation are
$$124,125,145,234,348,458,
568,256,236,138,128,278,678,237,137,167,136.$$ It is known that the
simplicial complex corresponding to this triangulation  is not
shellable (not even constructible), but it is Cohen-Macaulay, see
\cite{Ha}, and it has the following partition:
$$[\emptyset, 124 ], [3, 234], [5, 145], [6, 236], [7, 137], [8, 348], [13,
138], [16, 136], [18, 128], $$ $$[25, 125], [27, 237], [28, 278],
[56, 256], [67, 167], [68, 568], [78, 678], [58, 458].$$ \noindent
Therefore we have again
$\mathrm{depth}(\Delta)=\mathrm{dim}(\Delta)=\mathrm{sdepth}(\Delta)=3$.
}
\end{Example}

\begin{Example}{\em  (The Cylinder) The ideal $I=(x_1x_4,x_2x_5,x_3x_6,x_1x_3x_5, x_2x_4x_6)\subset K[x_1,\ldots,x_6]$
is the Stanley-Reisner ideal of the  triangulation of the cylinder
shown in the next figure. The corresponding simplicial complex
$\Delta$ is Buchsbaum but not Cohen-Macaulay.

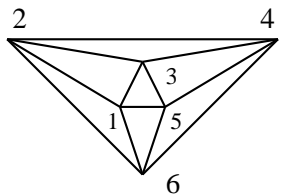
\psset{unit=0.3cm}
\begin{pspicture}(-11,6)(-11,17)

\put(9.2,11){\footnotesize 1} \put(12,11){\footnotesize 5}
\put(11.8,13){\footnotesize 3} \put(16,15.5){$4$} \put(5,15.5){$2$}
\put(11.8,8.2){$6$}

\psline(4.77,15)(16.77,15)\psline(4.77,15)(10.77,9)\psline(16.77,15)(10.77,9)

\psline(10.77,14)(9.77,12)\psline(9.77,12)(11.77,12)\psline(11.77,12)(10.77,14)

\psline(4.77,15)(10.77,14)\psline(4.77,15)(9.77,12)

\psline(16.77,15)(10.77,14)\psline(16.77,15)(11.77,12)

\psline(10.77,9)(9.77,12)\psline(10.77,9)(11.77,12)

\end{pspicture}

\noindent The  facets of $\Delta$  are $123, 126, 156, 234, 345,
456$, and it has the following partition:
$$[\emptyset, 123], [4, 234], [5, 345], [6, 456], [15, 156], [16, 126], [26,26].$$

\noindent Therefore we have
$\mathrm{depth}(\Delta)=\mathrm{sdepth}(\Delta)=2<3=\mathrm{dim}(\Delta)$.
Although $\Delta$ is not partitionable,  $I_\Delta$ is a Stanley
ideal.

}
\end{Example}

\section*{Acknowledgments}

This paper was prepared during the third author's visit of the
Universit\"at Duisburg-Essen, where he was on sabbatical leave from
the University of Tehran. He would like to thank Deutscher
Akademischer Austausch Dienst (DAAD) for the partially support. He
also thanks the authorities of the Universit\"at Duisburg-Essen for
their hospitality during his stay there.

\end{document}